\documentclass[preprint,noinfoline]{imsart}

\usepackage{amsthm,amsmath,amssymb,natbib}

\newtheorem{theorem}{Theorem}
\startlocaldefs
\def\bb#1{\mathbb{#1}}
\def\bbp{\bb{P}}
\def\bbf{\bb{F}}
\def\bbn{\bb{N}}
\def\en{\infty}
\def\summ#1#2#3{\sum_{#1=#2}^{#3}}
\def\mby{\boldsymbol{y}}
\def\cov{\mathop{\rm cov}}

\endlocaldefs

\begin{document}

\baselineskip=22pt

\begin{frontmatter}

\title{The Poisson Compound Decision Problem Revisited}
\runtitle{Poisson Compound}


\author{\fnms{Lawrence D.} \snm{Brown,}\thanksref{t1}\ead[label=e1]{lbrown@wharton.upenn.edu}}
\address{University of Pennsylvania; \printead{e1}}
\thankstext{t1}{Research supported by an NSF grant.}
\affiliation{University of Pennsylvania}
\author{\fnms{Eitan} \snm{Greenshtein}\ead[label=e2]{eitan.greenshtein@gmail.com}}
\address{Israel Census Bureau of Statistics; \printead{e2}}
\affiliation{Israel Census Bureau of Statistics}
\and
\author{\fnms{Ya'acov} \snm{Ritov}\thanksref{t3}\ead[label=e3]{yaacov.ritov@gmail.com}}
\address{The Hebrew University of Jerusalem; \printead{e3}}
\thankstext{t3}{Research supported by an ISF grant.}
\affiliation{The Hebrew University of Jerusalem}

\runauthor{Brown, Greenshtein, Ritov}

\begin{abstract}
The compound decision problem for a vector
of independent Poisson random variables with possibly different means has a half-century old solution. However, it appears that the classical solution needs smoothing adjustment. We discuss three such adjustments. We also present another approach that first transforms the problem into the normal compound decision problem. A simulation study shows the effectiveness of the procedures in improving the performance over that of the classical procedure.
A real data example is also provided. The procedures depend on a smoothness parameter, that can be selected using a
non-standard cross-validation step which is of  independent interest.  Finally, we mention some asymptotic results.
\end{abstract}



\end{frontmatter}

\section{Introduction}
\label{sec:introduction}

In this paper we consider the problem of estimating a vector
$\boldsymbol{\lambda}=(\lambda_1,\dots,\lambda_n)$, based on
observations $Y_1,\dots,Y_n$, where $Y_i \sim Po(\lambda_i)$
are independent. The performance of an estimator
$\hat{\boldsymbol{\lambda}}$ is evaluated based on the risk
\begin{equation} E_{\boldsymbol{\lambda} } ||\hat{
\boldsymbol{\lambda} }-\boldsymbol{\lambda}||^2,
\label{eqn:risk1}  \end{equation}
which corresponds to the loss function
$$ L_2(\boldsymbol{\lambda},\hat{ \boldsymbol{\lambda}})
= \sum ( {\lambda_i}- \hat{\lambda}_i)^2.$$

Empirical Bayes (EB) is a general approach to handle compound
decision problems. It was suggested by Robbins, see (1951,
1955);
see Copas (1969) and Zhang (2003) for review papers.
The improvement that empirical Bayes methods yield over more classical,
e.g, mle, methods is especially prominent in inference for high dimensional data. Thus the empirical
Bayes method has become especially relevant in recent years; see e.g., the enthusiastic advocation
for Empirical Bayes usage and relevance in Efron (2003).

Assume that
$\lambda_i$, $i=1,\dots,n$  are realizations of i.i.d.
$\Lambda_i$, $i=1,\dots,n$, where $\Lambda_i \sim G$.
Then a  natural approach is to use the Bayes procedure:
\begin{equation}
\delta^G= \mathop{\rm argmin}_\delta E_G(\delta(Y) -
\Lambda)^2,
\label{eqn:risk}
\end{equation}
and estimate $\boldsymbol{\lambda}$ by
$\hat{\boldsymbol{\lambda}}=(\delta^G(Y_1),\dots,\delta^G(Y_n))$.
When $G$ is completely unknown, but it is assumed that
$\lambda_1,\dots,\lambda_n$ are i.i.d., then it may be possible
to estimate $\delta^G$  from the data $Y_1,\dots,Y_n$, and
replace it by some $\hat\delta^G$.

 Optimal frequentist
properties of  $\delta^G$ in the context of the compound
decision problem,
are described in terms of optimality within the class of simple
symmetric decision functions. See the recent paper by Brown and
Greenshtein (2009) for a review of the topic. The optimality of
empirical Bayes decision rules  within the larger  class of
permutational invariant decision functions is shown in
Greenshtein and Ritov (2009).

The Bayes procedure $\delta^G$ has an especially simple form in
the Poisson setup. In this case there is also a simple and
straightforward estimator $\hat{\delta}^G$ for $\delta^G$.
Denote by  $P$ the joint distribution of $(\Lambda,Y)$, which
is induced by $G$. The Bayes estimator of $\lambda_i$ given an
observation $Y_i=y$, is:

\begin{equation}\begin{split}
\delta^G(y) \equiv E(\Lambda_i|Y_i=y)&= \frac{\int \lambda
P(Y_i=y|\Lambda_i=\lambda)dG(\lambda)}{ \int
P(Y_i=y|\Lambda_i=\lambda)dG(\lambda)}
\\
&= \frac{ (y+1)P_Y(y+1)}{P_Y(y)},
\label{eqn:Bayes}
\end{split}\end{equation}
where $P_Y$ is the marginal distribution of $Y$ under $P$.
Given $Y_1,\dots,Y_n$, we may estimate $P_Y(y)$ trivially by
the empirical distribution: $ \hat{P}_Y (y) = { \#  \{ i| Y_i=y
\} } / {n}.$
We obtain the following Empirical Bayes  procedure
\begin{equation}
\hat{\delta}^G(y) = \frac{ (y+1)\hat{P}_Y(y+1)}{\hat{P}_Y(y)}.
\label{eqn:EB}
\end{equation}

This estimator was originally proposed in Robbins (1955). It is still
the ``default''/``classical'' empirical Bayes
estimator
in the Poisson situation.
Various theoretical results established in the above-mentioned
papers and many other papers imply that
as $n \rightarrow \infty$, the above procedure  will have
various optimal properties.
This is very plausible, since as $n \rightarrow \infty$,
$\hat{P}_Y \rightarrow P_Y$
and thus $\hat{\delta}^G \rightarrow \delta^G$. However, the
convergence may be very slow, even in
common situations as demonstrated in the following example, and one might want to improve the above $\hat{\delta}^G$. This is the  main purpose of this work.

{\bf Example 1:} Consider the case where $n=500$ and
$\lambda_i=10$, $i=1,\dots,500$. The Bayes risk of $\delta^G$
for a distribution/prior $G$ with
all its mass concentrated at 10 is, of course, 0. The risk of
the naive procedure which estimates $\lambda_i$ by $Y_i$,
equals the sum of the variances, that is, $10 \times 500=5000$.
In 100 simulations
we obtained an average loss of 4335 for the procedure
(\ref{eqn:EB}), which  is
not a compelling improvement over the naive procedure, and is very
far from the Bayes risk.

We will improve  $\hat{\delta}^G$ mainly through ``smoothing''.
A non-trivial improvement is also obtained by imposing monotonicity
on the estimated decision function. By imposing monotonicity
without any further smoothing step, the average total loss in the above
example in 100 simulations is reduced to  301. By implementing
the procedure of Section 2 with a
suitable smoothing parameter ($h=3$) and imposing monotonicity the
average loss is reduced further to 30.
Early attempts to improve (\ref{eqn:EB})
through smoothing,
including imposing monotonicity, may be found in Maritz (1969)
and references there, see also Park (2011)
for further references and for an interesting application.

The rest of the paper is organized as follows.
In Section \ref{sec:adjusting} we will suggest  adjustments and
improvements of  $\hat{\delta}^G$.
In Section \ref{sec:normal} we describe the alternative
approach of
transforming the Poisson EB problem to a normal EB problem,
using a variance stabilizing transformation.
In Section \ref{sec:loss} we discuss some decision-theoretic
background, and in particular we examine loss functions other
than  squared-error loss. In Section \ref{sec:simulations}
we discuss the  above mentioned two approaches and compare them
in a simulation study. Both approaches involve a choice of a
``smoothing-parameter''.
For our new approach a choice based on cross-validation is suggested in Section
\ref{sec:XV}.
In Section \ref{sec:realData} we present an analysis of real
data describing frequency of car accidents.
In Section \ref{sec:NPMLE}
a further  approach which estimates
$\delta^G$ using a nonparametric MLE is discussed. Finally, in Section 9, we
study some asymptotic properties of the classical Robbins' estimator.

\section{ Adjusting the classical Poisson empirical Bayes
estimator }
\label{sec:adjusting}

Section 1 describes the Bayes decision function
$\delta^G$ and its straightforward estimator $\hat{\delta}^G$.
Surprisingly, it was found empirically (see, Example 1) that even for $n$
relatively large, when the empirical distribution is close to
its expectation, the estimated decision function should be
smoothed.
We discuss in this section how this estimator can be improved.
The improvement involves three steps,
which finally define an adjusted Robbins estimator.

\subsection{Step 1}

Recall the joint probability space defined on $(Y,\Lambda)$. We
introduce a r.v. $N \sim Po(h)$,
where $N$ is independent of $Y$ and $\Lambda$. Let $Z=Y+{ N}$.
Consider the suboptimal decision function
\begin{equation}\label{dh1}
\delta_{h,1}(z) \equiv E(\Lambda | Z=z) = E(\Lambda + h|
Z=z)-h.
\end{equation}
The above is the optimal decision rule, when obtaining
the corrupted observations $Z_i=Y_i+N_i$, $i=1,\dots,n$ instead of
the
observations $Y_1,\dots,Y_n$. The ``corruption parameter'' $h$
is a selected parameter,
 referred to as the ``smoothing parameter''.
In  general, we will select a
smaller $h$ as $n$ becomes larger. See Section \ref{sec:XV} for further discussion on the choice of $h.$
Motivated by \eqref{dh1} and reasoning similar to  \eqref{eqn:EB}, we define
$\hat{\delta}_{h,1}$ as:
\begin{equation}
\hat{\delta}_{h,1}(z)=  \frac{
(z+1)\tilde{P}_Z(z+1)}{\tilde{P}_Z(z)} - h,
\label{eqn:aadj1}
\end{equation}
when $\tilde{P}_Z(z) >0$; $\hat{\delta}_{h,1}(z)=0$ otherwise.

Here the distribution  $\tilde{P}_Z(z)$ is defined by
\begin{equation}\label{PZz}
\tilde{P}_Z(z)=  \sum_{i=0}^z \hat{P}_Y(i) \times \exp(-h)
\frac{h^{z-i}}{(z-i)!}.
\end{equation}
Note that $\tilde P_Z( z)$
as defined in \eqref{PZz}
involves observation of $Y$ through the quantity $\hat P_Y( y)$ that appears inside its definition. It
is---in general---a much better estimate of $ P_Z( z)$ compared to the empirical distribution function $\#\{i:\;Z_i = z\}$.

\subsection{ Step 2.}

There is room for considerable improvement of $\delta_{h,1}$. Note that $\delta_{h,1}$  is applied to the
randomized observation $Z_i$.  Therefore, the natural next
adjustment is Rao-Blackwellization of the estimator. Define
\begin{equation} \hat{\delta}_{h,2}(y)=E_{\cal N}
(\hat{\delta}_{h,1}(y+ {\cal N})), \label{eqn:adj2}
\end{equation}
for ${\cal N} \sim Po(h)$, which is independent of the
observations $Y_i, \; i=1,\dots,n$. That is,
\begin{equation*}
 \hat\delta_{h,2}(y) = e^{-h}\sum_{j=0}^{\infty} \frac{h^j}{j!}
\hat{\delta}_{h,1}(y+j).
\end{equation*}
Note that for a given $y$, the value of $\hat\delta_{h,2}(y)$
depends on all of $\hat P_Y(0),\hat P_Y(1),\dots$, although
mainly on the values in the nearby neighborhood of $y$.

\subsection{Step 3}

Finally after applying adjustments 1 and 2 we obtain a decision
function which is not necessarily
monotone. However, because of the monotone likelihood ratio
property of the Poisson model, $\delta^G$ is monotone.
A final adjustment is to impose monotonicity on  the decision
function $\hat{\delta}_{h,2}$. We do it through applying
isotonic regression by the pooling adjacent violators, cf.
Robertson, Wright, and Dykstra (1988).
Note, the monotonicity is imposed  on $\hat{\delta}_{h,2}$
confined to the domain $D(Y)\equiv \{y:\; Y_i=y \text{ for some }i=1,\dots,n\}$. To be more explicit, an estimator is isotonic if
 \begin{equation}
 \label{iso}
 y_i,y_j \in D(Y) \text{ and } y_i\leq y_j \Rightarrow \delta(y_i)\leq \delta(y_j),
 \end{equation}
and $\delta_{h,3}$ is isotonic and satisfies
$$\sum_{i=1}^n\bigl(\hat\delta_{h,3}(y_i) - \hat\delta_{h,2}(y_i)\bigr)^2 = \min\Bigl\{\sum_{i=1}^n\bigl(\hat\delta(y_i) - \hat\delta_{h,2}(y_i)\bigr)^2 : \delta \text{ satisfies  \eqref{iso}}\Bigr\}.$$

We  obtain the final decision function $\hat{\delta}_{h,3}$,
after this third step.

In order to simplify notations we denote:
${\Delta}_h \equiv \hat{\delta}_{h,3} .$ This is our adjusted
Robbins estimator.

\subsection{Discussion}

We now further discuss the above approach. Step 1 of this approach transforms the original
problem of
estimating the decision function in the Bayesian problem where $\Lambda \sim G$, to an
auxiliary problem
of estimating the decision function in a problem with $\Lambda' \sim G'$, where
$G'(\lambda+h)=G(\lambda)$. The estimation of the decision function in the auxiliary problem
is done through an adaptation of Robbins' classical estimator. Indeed, note that $(\hat{\delta}_{h,1}+h)$
is an estimator of the Bayes procedure in the auxiliary problem, using
(an adapted) Robbins' method.

Let $B(G)$ and $B(G')$ be the Bayes risk in the original and in the auxiliary problem.
Obviously $B(G') \geq B(G)$ since the original experiment dominates the auxiliary one.
Furthermore, as precisely argued in the final section, in both  the original
and the auxiliary problems the
difference between the average risk per coordinate of Robbins' procedure
and the Bayes risk is of order $o(\log(n)^2/n)$.
Hence, the average risk per coordinate, of our final procedure $\Delta_h$, is bounded below by $B(G)$ and bounded above by $B(G')+o(\log(n)^2/n)$.
For a fixed $h$  in non-trivial situations $\delta^G - \delta^{G'}$ does not converge to zero, and thus
for large enough $n$ our adjusted Robbins' procedure performs worse than the original Robbins'
procedure. However, our simulations show that the asymptotics might `kick-in' only for a very large $n$, and adjusting Robbins' procedure may be very helpful even for large values of $n$.
Estimating the decision function in the auxiliary problem could be much more efficient,
compared to estimating the decision function in the original problem,
even for large $n$.
The above heuristically implies, i) as $n$ grows we should use smaller $h$ ii) for distributions $G$ closer
to 0 (i.e., with smaller values) we might want to apply smaller values of $h$, since we expect a larger
difference between $B(G)$ and the upper bound $B(G')$.

The Rao-Blackwellization in Step 2
is especially needed when $h$ is not small. Note again, that  $(\hat{\delta}_{h,1}+h)$ is an estimator of the decision function $\delta^{G'}$, which might be very different than $\delta^G$ when $h$ is not small. In Step 2 we transform the original observations $Y_1,...,Y_n$, to
$Z_1,...,Z_n$ which are distributed according to the observations in the
auxiliary problem to which $\delta^{G'}$
corresponds, it is then averaged over all possible $Z_i, \; i=1,...,n$, in order to obtain a Rao-Blackwell improvement.

The choice of the PAV algorithm for the smoothing Step 3 is heuristically natural and convenient.
See, for example, Mammen (1991). But there could be other ways to carry out this step. See Koenker and Mizera (2012) for a recently proposed and interesting approach for monotonization and estimation. Our experience is that monotonization is particularly useful when $h$ is small since for larger $h$ the smoothing in the first two steps typically yields an estimator that is already very close to being monotone.

\bigskip

Finally we remark on a curious discontinuity property of
$\Delta_h$.
The function $\Delta_h$ is a random function, which depends on
the
realization $\boldsymbol{y}=(y_1,\dots,y_n)$. In order to
emphasize it we write here $\Delta_{\boldsymbol{y},h} \equiv
\Delta_h$.
Consider the collection  of functions parameterized by $h$,
denoted
$\{ \Delta_{\boldsymbol{y},h}(y) \}$. It is evident from the
definition of (\ref{eqn:aadj1}),
that $ \Delta_{\boldsymbol{y},h}(y)$ does not (necessarily)
converge to
$\Delta_{\boldsymbol{y},0}(y)$ as $h$ approaches 0, even for
$y$ in the range $y_1,\dots,y_n$. This will happen whenever
there is a gap in the range of $y$. Suppose, for simplicity
that $\hat P_Y(y)=0$, while $\hat P_Y(y-1),\hat P_Y(y+1)>0$.
Then,  $\lim_{h\to0}\hat\delta_{h,1}(y-1) = 0$, and
$\lim_{h\to0}h\hat\delta_{h,1}(y)=(y+1)\hat P_Y(y+1)/\hat
P_Y(y-1)$. Hence
\begin{eqnarray*}
 \lim_{h\to0} \hat\delta_{h,2}(y-1)&=&
\lim_{h\to0}E\Bigl(\hat\delta_{h,1}(y-1+N)
\bigl|y_1,\dots,y_n\Bigr)
 \\
 &=& \lim_{h\to0}
\bigl((1-h)\hat\delta_{h,1}(y-1)+h\hat\delta_{h,1}(y)\bigr)
 \\
 &=&  (y+1)\hat P_Y(y+1)/\hat P_Y(y-1),
\end{eqnarray*}
which is strictly different from $\hat{\delta}_{0,2}(y)=0$. Suppose, more generally, that $\hat P_y(y)>0$ and $\hat P_Y(y+j_0)>0$ for some $j_0>1$, but $\hat P(y+j)=0$ for $j=1,\dots,j_0-1$. Then one can check directly from the definition that $\lim_{h\to0}\hat \delta_{h,2}=(y+j_0)\hat P_Y(y+j_0)/\hat P_Y(y)$. Note that in such a situation $\hat \delta^G(y)=0$. Hence $\hat \delta_{h,2}(y)$ for small to moderate $h$ seems preferable to $\hat\delta^G(y)=\hat\delta_{0,2}(y)$ in such gap situations.

This phenomena is reflected in our simulations in Section
\ref{sec:simulations}, especially in Table \ref{Table5}.

Another curious feature of our estimator is when applied on $y_{max}=\max \{ Y_1,...,Y_n \}$.
It may be checked that: $\hat{\delta}_{h,2}(y_{max})= (y_{max}+1)h + O(h^2)$. When $h$ is small
so that $(y_{max}+1)h \ll  y_{max}$, this would introduce a significant bias. Hence, choosing very
small $h$, might be problematic, though this bias is partially corrected  through the isotonic regression.

\section{Transforming the data to normality.}
\label{sec:normal}
The emprical Bayes approach for the analogous normal problem
has also been studied for a long time.
See the recent papers of Brown and Greenshtein (2009) and of
Wenhua and  Zhang (2009)
and references there.
The Poisson problem and the derivation of (\ref{eqn:EB}) are
simpler and were  obtained by Robbins at a very
early stage, before the problem of density estimation, used in
the normal empirical Bayes procedure, was addressed.
In what follows we will describe the modification of
the normal method to the Poisson problem.

In the normal problem we observe $Z_i \sim N(M_i,\sigma^2)$,
$i=1,\dots,n$ where $M_1,\dots,M_n$ are i.i.d. random variables
sampled from $G$ and
the purpose is to estimate $\mu_1,\dots,\mu_n$ the realizations
of $M_1,\dots,M_n$. The application of the normal EB procedure
to the Poisson problem
has a few simple steps. First
transform the Poisson variables $Y_1,\dots,Y_n$  to the
variables $Z_i= 2*\sqrt{Y_i+q}$. Various recommenations
for q are given in the literature, the simplest and most common
one is $q=0$, but
the choice $q=0.25$ was recommended by
Brown et. al. (2005, 2009). Thus treat $Z_i$'s as
(approximate) normal variables with variance $\sigma^2=1$
and mean $2*\sqrt{\lambda_i}$, and estimate their means by
$\hat{\mu}_i$, by applying normal empirical Bayes
technique;
specifically, $\hat{\mu}_i= {\delta}_{N,h}(Z_i)$, as defined in
\eqref{eqn:ENEB} below. Finally estimate $\lambda_i=E Y_i$,
by  $\hat{\lambda}_i=\frac{1}{4}\hat{\mu}_i^2$.

We will follow the approach  of Brown and Greenshtein (2009).
Let $$ g(z)=  \int \frac{1}{\sigma}\varphi \Bigl( \frac{
z-\mu}{\sigma} \Bigr)dG(\mu).$$
It may be shown that the normal Bayes procedure denoted $\delta^G_N$, satisfies:
\begin{equation} \delta_N^G(z)= z + \sigma^2 \; \frac{g'(z)}{g(z)}.  \label{eqn:NEB} \end{equation}
The procedure studied in Greenshtein and Brown  (2009),
involves an estimation of $\delta_N^G$,
by replacing $g$ and $g'$ in \eqref{eqn:NEB} by their kernel
estimators  which are derived through
a normal kernel with bandwidth $h$. Denoting the kernel
estimates by $\hat{g}_h$ and $\hat{g}'_h$
we obtain the decision function, $(Z_1,\dots,Z_n)\times z
\mapsto R$:
\begin{equation} {\delta}_{N,h}(z)= z + \sigma^2 \;
\frac{\hat{g}_h'(z)}{\hat{g}_h(z)}.  \label{eqn:ENEB}
\end{equation}

One might expect this approach to work well in setups where
$\lambda_i$ are large, and hence,
the normal approximation to $Z_i=\sqrt{Y_i+q}$ is good. In
extensive simulations the above approach was
found to also work well  for configurations with moderate and
small values of $\lambda$.
In many cases it gave results comparable to the adjusted  Poisson EB
procedure.


\noindent\textbf{Remark}
In the paper of Brown and Greenshtein the estimator
$\delta_{N,h}$
as defined in (\ref{eqn:ENEB}) was studied.  However, just
as in the
Poisson case, it is natural to  impose monotonicity.
In the simulations of this paper we make this
adjustment using isotonic regression.
Again, the monotonicity is imposed on $\delta_{N,h}$
confined to the range $\{y_1,...,y_n \}$.
We denote the adjusted estimator by
\begin{equation*} \Delta_{N,h}.  \end{equation*}

\section{The  loss functions.}
\label{sec:loss}

The estimator $\delta_{N,h}(Z_i)=\hat{\mu}_i$, may be
interpreted as an approximation of the nonparametric EB
estimator for $\mu_i \equiv 2\sqrt{\lambda}_i$, based on
the (transformed) observations $Z_i$  under the loss
$L(\boldsymbol{\mu},{\boldsymbol{a}})=||{\boldsymbol{\mu}}-\boldsymbol{a}||^2$,
for the decision $\boldsymbol{a}=(a_1,\dots,a_n)$.
Thus, $\frac{1}{4}\hat{\mu}_i^2$ may be interpreted as the
approximation of the empirical Bayes estimator
for $\lambda_i$, under the
loss
$$ L_H(\boldsymbol{\lambda},{\boldsymbol{ a}})= \sum
(\sqrt{\lambda_i}- \sqrt{a_i} )^2=-2\log(1-D_H^2),$$
where $D_H$ is  to the  Hellinger
distance between the distributions $\prod Po(\lambda_i)$
and $\prod Po(a_i)$.

Some papers  that  discuss the problem of estimating a
vector of Poisson means are
Clevenson and Zidek (1975), Johnstone (1984), Johnstone and
Lalley (1984) and Fourdinier and Robert(1995).
Those and other works suggest that a particularly natural
loss function
in addition to $L_H$ and $L_2$, denoted $L_{KL}$ is
$$ L_{KL}(\boldsymbol{\lambda},\hat{\boldsymbol{
\lambda}})= \sum \frac{ ({\lambda_i}- \hat{\lambda}_i)^2
}{\lambda_i}.$$
Note, $L_{KL}$ also corresponds to the local
Kulback-Leibler distance between the distributions.

From an empirical Bayes perspective,
the optimal decisions
that correspond to those three loss functions may have more
and less similarity, depending on the configuration.
For example, when the prior $G$ is concentrated on a point
mass, the Bayes procedures corresponding to
those 3 loss functions are obviously the same. Since the
$L_{KL}$ loss is of a special importance, we will briefly
describe how our analysis can be modified to handle it.
As in the case of  $L_2$ loss, one may obtain that the Bayes
decision under the
$L_{KL}$ loss is given for $y\ge 1$ by:   $$ \frac{y
P_Y(y)}{P_Y(y-1)}.$$ The decision for $y=0$ denoted
$\hat{\lambda}(0)$,
is:
\begin{eqnarray*}
\hat{\lambda}(0) &=& \arg\min_a \int
\frac{(\lambda-a)^2}{\lambda}e^{-\lambda} dG(\lambda)
\\
&=& \frac{\int e^{-\lambda} dG(\lambda)}{\int \lambda^{-1}
e^{-\lambda} dG(\lambda)}.
\end{eqnarray*}
In particular, $\hat\lambda(0)=0$ if $G$ gives a positive
probability to any neighborhood of 0.

The decision for $y \geq 1$ may be estimated as in
(\ref{eqn:EB}) together with
the three adjustments suggested in Section
\ref{sec:adjusting}, along the same lines.
However, we still need to approximate the Bayes decision
$\hat{\lambda}(0)$. Note however, that if $G$ has a point mass at 0, however small, the risk will be infinite unless $\hat\lambda(0)=0$. This is the only safe decision, since we cannot ascertain that there is no mass at 0.

Note, defining $Z=Y+{\cal N}$, ${\cal N} \sim Po(h)$ under the KL loss
as in Step 1 in the squared loss, might introduce instability due to small values of
$\tilde{P}_Z(z-1)$ in the denominator of
$\tilde{P}_Z(z)/\tilde{P}_Z(z-1)$, e.g., for $z=\min\{Z_1,...,Z_n\}$. One might want to define the "corrupted" variable
alternatively, as $Z \sim B(Y,p)$. Then $Z \sim Po(p \lambda)$, when $Y \sim Po(\lambda)$.
Our smoothing/corrupting parameter is $p$. We skip the details of the analogs of steps 1-3.

Throughout the rest of the paper, we consider and evaluate
procedures explicitly  only under the $L_2$ loss.

\section{Simulations}
\label{sec:simulations}

In this section we provide some simulation results which approximate
the risk of various procedures as defined in (\ref{eqn:risk1}).
Specifically for various {\it fixed} vectors $\boldsymbol{\lambda}=(\lambda_1,...,\lambda_n)$, we estimate
$ E_{\boldsymbol{\lambda}} \sum ( \Delta_h(Y_i)-\lambda_i)^2$ and
$ E_{\boldsymbol{\lambda}} \sum ( \Delta_{N,h}(Y_i)-\lambda_i)^2$, for various values of $h$.
The results are reported in tables bellow, each entry in those tables is based on 1000 simulations. 

It is known that for fixed vector $\boldsymbol{\lambda}=(\lambda_1,...,\lambda_n)$
a good benchmark and a lower  bound for the risk of our suggested procedures is
$nB(\boldsymbol{\lambda})$; here $B(\boldsymbol{\lambda})$ is the Bayes risk for the problem where we observe
$\Lambda \sim G$, where $G$ is the empirical distribution which is defined by $\lambda_1,...,\lambda_n$.
See Greenshtein and Ritov (2009) for a general investigation and discussion of this relation.

As already seen in Example 1, adjusting the
classical
non parametric empirical Bayes estimator can yield a significant improvement
in the risk. Significant improvement also occurs in a range of parameter
configurations, as exemplified by those in the following tables.
The normal empirical Bayes method of Section 3 works nearly as well in many of those
configurations, but seems less suited to tightly clumped configurations like those in Tables 3 and 4.
We were somewhat surprised to find that this normal method does compare reasonably well even when there are some small values of $\lambda$ as in Table 2. Simulations for the normal method were performed with both $q=0$ and
$q=1/4$, as variance stabilizers.
In every case the results for $q=1/4$ were between $2\%$ and $5\%$ better than those for $q=0$. So, we report only on those with $q=1/4$.


We elaborate on Table \ref{Table1}. The reading of the other
tables is similar. In Table 1 we study risks of our procedure, $\Delta_h$, and of
of $\Delta_{N,h}$ for various values of $h$. The risks for this table are computed when
$\lambda_1,...,\lambda_{200}$ are equally spaced between 5 and 15. In practical settings the smoothing parameter, $h$, should be selected according to cross-validation or other method. In Section 6 we describe a new cross-validation method that seems to work well in the present context. In Table 1 the risks of $\Delta_h$ and $\Delta_{N,h}$ under the perspective best choices of $h$ are shown in bold-face. The second row of the table shows the risk of
$\hat{\delta}_{h,2}$. This procedure does not involve the isotonic monotonization step.
This is included for the purpose of comparison in order to show the beneficial effect of this final step of our procedure.

Note that $\hat{\delta}_{0,2}$ is the classic Robbins' procedure. Its risk is much larger than
is available from $\Delta_h$ or $\Delta_{N,h}$. The risk of $\Delta_0$ is is that of the classic procedure followed by the monotonization step and, as can be seen, this step considerably reduces the risk. However, as $h$ increases, the procedure $\hat{\delta}_{h,2}$
becomes more nearly monotone and as can be seen from the table the monotonization step becomes less important in decreasing the risk.

For purposes of comparison we note that the risk of the naive procedure is 1500 and the risk of the Bayes procedure for the setting of the table is approximately 880.

Our simulations were done using R (2008) software;
monotonicity is imposed on all the estimators, as described in Step 3,
through the `isoreg' R-procedure.

An observed  advantage, of the adjusted Robbins' method over the transformed normal method, is its stability with respect to the
chosen smoothing parameter $h$. This appears in Table 1 and is even more apparent in some of the
subsequent tables.

\begin{table}
\caption{Different EB procedures for
$\lambda_1,\dots,\lambda_{200}$ that are
evenly spaced between 5 and 15}\label{Table1}
\vspace{1ex}\begin{center}
\begin{tabular}{|l|rrrrrrr|}
\hline&&&&&&&\\
 & h &        0
& 0.2
& 0.4
& 0.8
&1.8
&3 \\
$\Delta_h$&risk & 1114 & 1049 & { 1017} &  994 & { 965} & {\bf 958}
\\ $\hat{\delta}_{h,2}$  & risk &6714& 2656& 1623& 1162& 994& 964\\&&&&&&&\\\hline&&&&&&&\\
$\Delta_{N,h}$ & h    &    0.2
  & 0.3
  & 0.5
  & 0.7
  & 0.9
  & 1.2 \\
&risk & 1230 & 1099 & { 1013} & { \bf 997}  & 1046
& 1138
 \\&&&&&&&\\\hline

\end{tabular}
\end{center}
\end{table}

The model studied in Table \ref{Table2}   is of $\lambda_i$,
$i=1,\dots,200$ evenly spaced between $0$ and $5$. Comparing
the two halves of the table, one may see how well the normal
modification works even for such small value of  $\lambda_i$.

\begin{table}
\caption{Different EB procedures for
$\lambda_1,\dots,\lambda_{200}$ that are
evenly spaced between 0 and 5} \label{Table2}
\vspace{1ex}\begin{center}

\begin{tabular}{|l|rrrrrrr|}
\hline&&&&&&&\\
 &h &    0
& 0.5
& 1
& 1.8
&2.4
&3 \\
$\Delta_h$&risk & 248 & {\bf 229} &{ 232} &  242 & { 249} & 258
\\$\hat{\delta}_{h,2}$&risk&556& 305& 233& 243& 250& 259\\&&&&&&&\\\hline&&&&&&&\\
$\Delta_{N,h}$& h &  0.2
& 0.3
& 0.5
&0.8
& 1.0
& 1.4 \\
 &risk & 308 & 267 & 245 & {\bf 242 } & 254 &
291\\
&&&&&&&\\\hline
 \end{tabular}

\end{center}
\end{table}

Next, in Table \ref{Table3},  we study the case where
$\lambda_1=\dots=\lambda_{200}=10$.
 Here the advantage of the adjusted Poisson over the modified
normal is clear.

\begin{table}
\caption{Different EB procedures for
$\lambda_1=\dots=\lambda_{200}=10$.}\label{Table3}
\vspace{1ex}\begin{center}
 \begin{tabular}{|l|rrrrrrr|}
\hline&&&&&&&\\
&  h     &    0
   & 0.2
  & 0.4
& 1
&2
&3 \\
$\Delta_h$&risk & 253 & 121 & 90 &  54 & { 38} & {\bf 28}
\\$\hat{\delta}_{h,2}$&risk&3904&1215&570& 160& 72& 47\\&&&&&&&\\\hline&&&&&&&\\
$\Delta_{N,h}$ & h    &    0.2
   & 0.3
  & 0.5
  & 0.7
  & 0.9
 & 1.3 \\
 &risk & 330 & 197 & {\bf 180} &  265  & 442 &
808
\\&&&&&&&\\\hline
 \end{tabular}
\end{center}
\end{table}

Next we study the following situation where we have a few
large $\lambda_i$ values:
$\lambda_1=\dots=\lambda_{200}=5$, while
$\lambda_{201}=\dots=\lambda_{220}=15$. There is still a clear
advantage
of the adjusted Poisson over the modified normal. See Table
\ref{Table4}.  It seems that in this situation  the advantage of the modified Robbins procedure  over the  normal is due to the
poor tail approximation of the latter.

\begin{table}
\caption{Different EB procedures for
$\lambda_1=\dots=\lambda_{200}=5$, while
$\lambda_{201}=\dots=\lambda_{220}=15$.}\label{Table4}
\vspace{1ex}\begin{center}

\begin{tabular}{|lr|rrrrrr|}
\hline&&&&&&&\\
& h      &    0
   & 0.2
& 0.4
& 1.2
&2.0
 &3 \\
$\Delta_h$ &risk & 665 & 476 & { 471 }& {\bf 449} & { 462} & 483
\\$\hat{\delta}_{h,2}$&risk&10382& 3488& 1761&
720& 623& 599\\&&&&&&&\\\hline&&&&&&&\\
$\Delta_{N,h}$ & h   &    0.2
  & 0.3
 & 0.5
 & 0.9
 & 1.1
 & 1.4 \\
 &risk & 819 & 613 & {\bf 550} & 653 & 732 &
823\\
&&&&&&&\\\hline \end{tabular}

\end{center}
\end{table}

Finally we investigate a configuration with only $n=30$
observations spread
over a larger interval. The $\lambda_i$ are evenly spread
between 0 and 20.
For this configuration there is a slight advantage of the
modified normal procedure.
In order to demonstrate the discontinuity  of $\Delta_h$
mentioned in Remark 1, we approximated the risk
of $\Delta_h$ for $h=0.01$, based on 1000 simulations. The
approximated risk is 244, compared
 to 867, for $h=0$, this is also the minimal approximated risk
from the values of $h$ that we tried
 in Table \ref{Table5}.

\begin{table}
\caption{Different EB procedures for
$\lambda_1,\dots,\lambda_{30}$ that are evenly spread between 0
and 20.}\label{Table5}
\vspace{1ex}\begin{center}

\begin{tabular}{|l|rrrrrrr|}
\hline&&&&&&&\\
 &
h      &    0
 & 0.2
 & 0.4
& 1.2
 &2.0
 &3 \\
$\Delta_h$&risk & 867 & 256 & {\bf 249 }& { 256} & { 262} & 260
\\$\hat{\delta}_{h,2}$&risk&3190& 1452& 924& 384& 320& 281\\&&&&&&&\\\hline&&&&&&&\\
$\Delta_{N,h}$ &
 h   &    0.2
  & 0.3
 & 0.5
 & 0.9
  & 1.2
 & 1.4 \\
 &risk & 316 & 302 & 280 & 243 & {\bf 236} & {
239}
\\&&&&&&&\\\hline
 \end{tabular}

\end{center}
\end{table}

Finally, the standard error of the estimated risk in the range of smoothing parameters $h$,
is about 3 in Experiment 1, about 1 in experiments 2-4, and about 2.5 in Experiment 5.


 \section{Choosing  the smoothing-parameter by
Cross-validation}
\label{sec:XV}

 In this section we  suggest a
non-standard
cross validation method, and study its performance.  This method is explained in the Poisson context,
and then
 in the normal context. The same general idea works for other
cases where an observation
 may be factorized, e.g., for infinitely divisible
experiments.
 About factorization of experiments,
 see Greenshtein (1996) and references there.


Let $p\in(0,1)$, $p\approx 1$, and  let $U_1,\dots,U_n$ be independent given $Y_1,\dots,Y_n$, where $U_i\sim B(Y_i,p)$, $i=1,\dots,n$, are binomial variables. As is well known, one of the features of the Poisson distribution is that $U_i\sim Po(p\lambda_i)$, and $V_i\equiv Y_i-U_i\sim Po((1-p)\lambda_i)$, and they are independent given $\lambda_1,\dots,\lambda_n$. We will use the main sub-sample $U_1,\dots,U_n$ for the construction of the family of estimators (parameterized by $h$), while the auxiliary sub sample $V_1,\dots,V_n$ is used for validation.
The choice $p \approx 1$ is in order that the distribution of $U_i$ will be close to that of $Y_i$, $i=1,...,n$, thus estimation based on $U_i$ is similar to estimation based on $Y_i$.
Let $\hat\delta^*_h(\cdot)$, $h\in H$ be a family of estimators, based on $U_1,\dots,U_n$ such that $\hat\delta^*_h(U_i)$ estimates $p\lambda_i$, $i=1,\dots,n$.  Consider:
\begin{equation}
\begin{split}
    &\hspace{-.2em}
    \rho(h;\boldsymbol{U},\boldsymbol{V})
    \\&= \frac1n \sum_{i=1}^n \Bigl( \hat\delta^*_h(U_i) - p(1-p)^{-1}V_i\Bigr)^2
\\
    &=  \frac1n \sum_{i=1}^n \Bigl( \bigl( \hat\delta^*_h(U_i) - p\lambda_i\bigr) - p(1-p)^{-1}\bigl(V_i-(1-p)\lambda_i\bigr)\Bigr)^2
    \\
    &= \frac1n \sum_{i=1}^n  \bigl( \hat\delta^*_h(U_i) - p\lambda_i\bigr)^2 + R_n(h)+A_n
      ,
\end{split}
\end{equation}
where $A_n$ is a random quantity that does not depend on $h$, and has no importance to the selection of $h$, while
\begin{equation}
\begin{split}
    R_n(h)&= \frac {2p}{(1-p)n}   \sum_{i=1}^n  \bigl( \hat\delta^*_h(U_i) - p\lambda_i\bigr)\bigl(V_i-(1-p)\lambda_i\bigr).
\end{split}
\end{equation}
Since $V_1,\dots,V_n$ are independent and independent of $U_1,\dots,U_n$ given $\lambda_1,\dots,\lambda_n$:
\begin{equation}
\begin{split}
    E(R_n^2(h)|\boldsymbol{U},\boldsymbol{\lambda})&= \frac {4p^2}{(1-p)n^2}   \sum_{i=1}^n  \bigl( \hat\delta^*_h(U_i) - p\lambda_i\bigr)^2\lambda_i.
\end{split}
\end{equation}

We conclude that if  $(1-p)n/\max\{\lambda_i\}|H|\to\infty$, then
\begin{equation}
 \rho(h;\boldsymbol{U},\boldsymbol{V})
 = L(\hat\delta^*_h,p\boldsymbol{\lambda})+o_p(1),
  \end{equation}
uniformly in $h\in H$. Recall that the decision function $\hat\delta^*_h$ used in the above result, is the non-parametric empirical Bayes procedure based on $U_1,\dots,U_n$ and $\hat\delta^*_h(U_i)$ is estimating $p\lambda_i$. If also $p\to 1$, we suggest to use the value $h$ that minimizes $\rho(h;\boldsymbol{U},\boldsymbol{V})$, to construct a similar estimator based on the original sample $Y_1,\dots,Y_n$, estimating $\lambda_1,\dots,\lambda_n$.

$\rho(h;\boldsymbol{U},\boldsymbol{V})$, given the sample $Y_1,\dots,Y_n$ is a randomized estimator of the loss function. Once again we suggest in this paper to replace a randomized estimator by its expectation given the sample $E\Bigl( \rho(h;\boldsymbol{U},\boldsymbol{V})\Bigl|\boldsymbol{Y}\Bigr)$. This expectation can be estimated by a Monte Carlo integration---sampling $K$ i.i.d. samples of $\boldsymbol{U}$ and $\boldsymbol{V}$.

 For the normal model,
 $Z_i \sim N(\mu_i,1)$, $i=1,\dots,n$, let   $\epsilon_{i}
\sim N(0,1)$ be  auxiliary i.i.d. variables,  independent of
$Y_1,\dots,Y_n$. Define $U_{i}=Y_i +\alpha \epsilon_{i}$,
$V_{i}= Y_i - (1/\alpha)\epsilon_{i}$.
 Then $U_{i}$ and $V_{i}$ are independent both with mean
$\mu_i$, and with variances $1+\alpha^2$ and $1+
(1/\alpha^2)$
 correspondingly.  Again, $\boldsymbol{U}$ may be used for estimation
and $\boldsymbol{V}$ for validation, where $\alpha>0$, $\alpha\to0$.

\subsection{Numerical Study.}

\noindent{\bf Example 2:} Consider the configuration
$\lambda_1=\dots=\lambda_{200}=10$,  simulated in Table
\ref{Table3} Section \ref{sec:simulations}. In that table $h=3$ is
recommended with a noticeable advantage over  $h \leq 0.4$. We
applied the above cross validation procedure with $p=0.9$ on a
single
realization of $Y_i, \;i=1,\dots,200$.  We repeated the
cross-validation
process $K=10000$ times on this single realization for the
values $h\in\{0, 0.5, 1, 1.5, 2, 2.5, 3 \}$. The corresponding
numbers $\rho(h,\bf{U},\bf{V})$ (scaled by $(1-p)^2$)  were: 165.834, 164.862, 164.736, 164.457,
164.421, 164.286, 164.340.  Note that, the last numbers
represent mainly the variance of our validation variable, but
the success of the corresponding estimator
is also a factor. The numbers indicate that the choices $h=0,
0.5, 1$ are inferior, the formal recommended
choice is $h=2.5$, the second best is $h=3$.

We repeated the simulation on another single realization, again
$K=10000$, this time we took $p=0.85$. The corresponding
numbers
are: 220.562, 217.986, 217.706, 217.374, 217.209, 217.272, 217.247.
Again, the numbers indicate that the choices
$h=0, 0.5, 1$ are inferior. The formal recommended choice is
$h=2$, the second best is again $h=3$.

\bigskip

Finally, we extended the five experiments, studied in the previous section.
For each experiment we repeated 100 times the following simulation.
We took the six values of $h$ which are reported in the corresponding table
in  Section \ref{sec:simulations}, and in each of the 100 runs we chose the smoothing parameter
among the six candidates through implementing  the above cross validation method with $K=10000$
and $p=0.9$. Hence different values of $h$ were used for different realizations.
The  results we obtained for experiments 1-5 are  correspondingly:  944, 246, 30, 453, 258.
The simulated risks that correspond to the best individual smoothing parameter in each
experiment are: 958, 229, 28, 449, 249. The performance of the CV is quite impressive.

Note that  in Experiment 1 the simulated risk of the CV is actually smaller than all the risks that correspond to the individual six smoothing parameters. This improvement could be an artifact
of the simulation and not a real one, our simulations were too slow to make a confident statement.
However, such an improvement could be real since the CV method might choose a different `more suitable' smoothing parameter depending on the realization.

 \section{Real Data Example.}
\label{sec:realData}

 In the following we study an example based on real data about
car accidents with injuries
 in 109 towns in Israel in
 July 2008. The 109 towns are those that had at least one
accident with injuries in that period of time;
 in the following we ignore this selection bias.
 There were 5 Tuesdays, Wednesdays and Thursdays, in that
month. For Town $i$, let $Y_i$ be the total
 number of accidents with injuries in those 5 Wednesdays.
 Similarly, for Town $i$, let $Z_i$ be half of the number of
accidents
 with injuries in the corresponding
 Tuesdays and Thursdays. We modelled  $Y_i$ as independently
distributed $Po(\lambda_i)$.
 In the following we will report on the performance  of our
empirical-Bayes estimator for various smoothing
 parameters $h$. It is evaluated through the predictive squared error
$$ \hat{R}=\sum (Z_i-\Delta_h(Y_i))^2.$$
 The towns Tel-Aviv and Jerusalem had a heavy impact on the
risk and thus we excluded them from the analysis.
 The remaining data seems to have relatively low values of
$\lambda_i$, a case where the classical
 Poisson-EB procedure is expected to perform well, and indeed
it does. The range of $Y_i$ is 0-14, while
 $\sum Y_i=135$,
 and $\sum Y_i^2=805$.
 In this example, the classical Poisson-EB adjusted for
monotonicity (i.e., $h=0$),
 gave the best result. Applying a smoothing parameter $h>0$ is
slightly inferior
 based on the above empirical risk. Yet, it is re-assuring to
see how stable is the performance of $\Delta_h$, as
 $h$ varies.   The empirical loss for the naive procedure
estimating $\lambda_i$ by $Y_i$, is 240.
The   modified normal estimators
with $q=\frac{1}{4}$ and various values of $h$ was applied to
the data as well. Again a clear advantage of our class of
adjusted
Poisson procedures over the class of modified normal procedures
was observed. In particular, the former class is much more
stable  with respect to the choice of the smoothing parameter
$h$. The results are summarized in Table \ref{Table:traffic}.

 \begin{table}
 \caption{EB applied to traffic accident by
city}\label{Table:traffic}
 \begin{center}
\begin{tabular}{|l|rrrrrrr|}
\hline&&&&&&&\\
$\Delta_h$ &
 h   &    0
  & 0.5
  & 1
  & 1.5
  & 2
 & 3 \\
&
 $\hat{R}$ & {\bf 140} & 163 &  172 &  168 & 166 & 159  \\
 &&&&&&&\\\hline&&&&&&&\\
 $\Delta_{N,h}$ &h  &  0.2
& 0.6
 & 1

  & 2
  & 3
 & 4 \\
&
 $\hat{R}$  & 262 & 185 &  174 & {\bf 170}   & 183 & 202  \\
 &&&&&&&\\\hline
 \end{tabular}
\end{center}
\end{table}

\section{The nonparametric MLE}
\label{sec:NPMLE}

The nonparametric maximum-likelihood (NPMLE),
as suggested by Kiefer and Wolfowitz (1956), is an alternative
approach for estimating $\delta^G$. It yields, automatically, a
monotone and smooth decision function. See Jiang and Zhang
(2009) for the normal model. To simplify the discussion, we
will assume that $\lambda_1,\dots,\lambda_n$ are realizations of
 i.i.d. random
variables sampled from the distribution $G$. Obtaining
a NPMLE $\hat{G}$ for $G$, induces the estimator $\delta^{\hat{G}}$ for $\delta^G$.
We will refer  to $\delta^{\hat{G} }$ also as  ${\delta}_{KW}$.

Note that the NPMLE maximizes with respect to $G$,  the
likelihood function:

\begin{equation*}\begin{split}
 \frac1n \summ i1n \log p_G(y_i) &= \summ i0\en
\bbp_n(i)\log p_G(i)
 \\
 &=  \summ i0\en \bigl(\bar\bbf_n(i-1)-\bar\bbf_n(i)\bigr)
\log p_G(i)
 \\
 &= \log p_G(0) + \summ i0\en \bar\bbf_n(i) \log \frac{
p_G(i+1)}{p_G(i)}
 \\
 &= \log p_G(0) + \summ i0\en
\bar\bbf_n(i)\log\delta^G(i) + C(\mby).
 \end{split}\end{equation*}
where $\bbp_n$ is the empirical process,
$\bbp_n(i)=\bbp_n(\{i\})$, and $\bar\bbf_n(i)=\summ j{i+1}\en
\bbp_n(j)$ ($\bar\bbf_n(-1)=1$). That is, the likelihood function can be written as a direct function of the Bayes procedure.

Suppose $G$ is supported on $[a,b]$.
Extend
$$\delta^G(y)=\frac{\int \lambda^{y+1}e^{-\lambda}dG(\lambda)
}{ \int \lambda^ye^{-\lambda}dG(\lambda)},\quad y\in R_+.$$
Then, clearly, $\delta^G(y)\in[a,b]$. Moreover, it is monotone
non-decreasing with derivative
$\delta^{G'}(y)=\cov(\lambda,\log\lambda)\in [0,b\log b-a\log
a]$ (where the covariance is with respect to measure  $\lambda^y
e^{-\lambda}dG(\lambda)$ normalized)

It is well known that the NPMLE is discrete with point mass $g_1,\dots,g_k$ on $\lambda_1,\dots,\lambda_k$ say. It is easy to see that it satisfies
 \begin{equation*}
    \summ i1n \frac{\lambda_j^{y_i}}{y_i!p_G(y_i)}=e^{\lambda_j}\quad, j=1,\dots, k.
    \end{equation*}
Since the left hand side is a polynomial in $\lambda$ of degree $\max y_i$, and a polynomial of degree $q$ in $\lambda$ can  be equal to $\exp\{\lambda\}$ only $q$ times, we conclude that $k<\max y_i$ (a more careful argument can reduce the bound on the support size). Hence, it is feasible to approximate algorithmically the NPMLE. Pursuing the asymptotic properties of the NPMLE estimator is beyond the scope of this paper. We should mention that as we argue in Section \ref{sec:theory}, Robbins' estimator is weak only when $G$ is  sparse and discrete, exactly where the NPMLE seems to excel.


Koenker and Mizera (2012) further developed this idea for the normal case. They approximated
$\delta_{KW}$ directly
(i.e., not through approximating $\hat{G}$ first), utilizing the monotonicity property/constraint of $\delta_{KW}$ to define
a corresponding convex optimization problem. Then, using interior point methods and available softwares they derived
algorithmically very efficient approximations of ${\delta}_{KW}$.

We are indebted to the AE for the
following Table 7 provided to us. In the first line of the table, the risk of the approximated ${\delta}_{KW}$
is given for the 5 simulated numerical experiments presented in our simulation section. Those results are  based on 1000 simulations for each example.
The  second line in the table gives the simulated risk of $\Delta_h$ for the best
value of $h$ among those reported in Tables 1-5, the third line gives the estimator obtained through cross-validation, as given and described in Section \ref{sec:XV}.

The performances of the methods are very similar. An advantage  of our suggested
procedure is that it is rather elementary and does not require more sophisticated  optimization methods and software. Also, as described in Section \ref{sec:loss} our method may be modified and specialized to deal with other loss functions.  It may also prove to be more adaptable for
generalizations involving additional covariates such as were studied in the normal case by
Jiang and Zhang (2010), Cohen, Greenshtein and Ritov (2012),
Koenker and Mizera (2012). We hope to study this issue in the future.

An advantage of  (the approximation of)
${\delta}_{KW}$ is that it does not involve a choice
of a smoothing parameter $h$, and does not require cross validation.

 \begin{table}
 \caption{ Comparison with Kiefer and Wolfowitz estimator}\label{ KW}
 \begin{center}
\begin{tabular}{|lrrrrr|}
\hline&&&&&\\ &  Exp1
&   Exp2
&   Exp3
&   Exp4
&   Exp5 \\\hline&&&&&\\
KW estimator&   {958} & 228 &  39 &  434 & 263  \\
Best-$h$ &    958  & 229 &  28 &  449 & 249  \\
CV selection&     944  & 246 &  30 &  453 & 258   \\
 &&&&&\\\hline
 \end{tabular}
\end{center}
\end{table}

\section{ Asymptotics  for Robbins' Estimator.}
\label{sec:theory}
In this section we will investigate theoretically the performance of Robbins'
method. It will be shown that in the  usual asymptotical EB setup, where we observe i.i.d. $\Lambda_1,...,\Lambda_n$, $\Lambda_i \sim G$, and $G$ is non-degenerate, Robbins' procedure $\hat{\delta}^G$ is very efficient. This is because that its risk is within $O((\log n/\log\log n)^2)$ of the risk of the Bayes procedure which is of order
$O(n)$. Note however, that if $G$ is degenerate the risk of the Bayes procedure is zero, and achieving a risk of order $(\log n/\log\log n )^2$ rather than a zero risk might not be considered a "success", in particular the ratio of the risks in that case is
infinity. More generally when the sequence $\lambda_1,...,\lambda_n$ of the realized $\Lambda_i \; i=1,...,n$
is very "sparse", in the sense that only a very small fraction of it does not equal to $\lambda_0$
, then Robbins' procedure,
whose risk will be shown to be larger than the Bayes risk by $\kappa(\log n/\log\log n) ^2$ for appropriate $\kappa>0$, might not be considered efficient. Note, we use the term sparse for
$\lambda_0$ which does not equal necessarily zero; in fact, the case $\lambda_0=0$ is excluded from the following theorem and from the
discussion, to avoid technical difficulties.

In order to formally study asymptotics for such sparse setups we will consider a triangular array
where at stage $k$,
$G=G^k$. Typically we consider  $G^k \rightarrow G_0$ weakly, where $G_0$ may be degenerate
at $\lambda_0$, where the support of $G^k$ is bounded uniformly in $k=1,2,..$.

For simplicity we assume further that the sample size $M$ is a Poisson random variable with mean $\nu=\nu^k$. Asymptotic results will hold as $\nu^k \to\infty$.  This assumption simplifies considerably the proof, and has little significance for the interpretation of the result.     Let $\bbn_\nu(y)=\#\{i:1\le i\le M,Y_i=y\}$, $y=0,1,\dots$. Note that they are independent under the Poisson sample size, $\bbn_\nu(y)\sim Po(\nu P(y))$, where $P(\cdot)$ denotes the marginal probabilities of $Y$.   A proof for a fixed sample size would involve the binomial distribution $B(P(y),n)$ for  $\bbn_n(y)$ and conditional on $\bbn_n(y)$, $\bbn_n(y+1)\sim B\bigl(n-\bbn_n(y),P(y+1)/(1-P(y)\bigr)$, but otherwise would be very similar, though more cumbersome. Let
\begin{equation*}
\begin{split}
\delta^{G^k}(y)&= (y+1) \frac{P(y+1)}{P(y)},\quad y=0,1,\dots
\\
{\hat\delta^{G^k}}(y) &= (y+1)\frac{\bbn_\nu^k(y+1)}{\bbn_\nu^k(y)},\quad y=0,1,\dots
\end{split}
\end{equation*}
be the Bayes procedure and its Robbins' approximation.

In the sequel we will occasionally drop the superscript $k$ for simplicity.

Let $r(G,\delta)$ be the total Bayes risk of the estimator $\delta$ when $\lambda_1,\dots,\lambda_M$ are (given $M$) simple random sample from $G$.

Our main result in this section is the following theorem.

\begin{theorem}

Suppose that $\lim\inf G^k\bigl((\lambda_1,\infty)\bigr)>0$ and $\lim\inf G^k\bigl([0,\lambda_2)\bigr)=1$ for some $0<\lambda_1<\lambda_2<\infty$ . Then   $(r(G^k, \hat{\delta}^{G^k}) - r(G^k,\delta^{G^k}))(\log\log \nu/\log \nu)^2 $ is bounded from above and away from 0.

\end{theorem}

\bigskip


\begin{proof}
 The risk of Robbins' procedure $\hat{\delta}^G$ is given by
\begin{equation*}
\begin{split}
r(G,\hat \delta^G)&= E\sum_{y=0}^\infty \bbn_\nu(y) E\Bigl(\bigl({\hat\delta^G}(y)- \Lambda\bigr)^2|Y=y\Bigr)
\\
&=  E\sum_{y=0}^\infty \bbn_\nu(y) \bigl({\hat\delta^G}(y)- \delta^G(y)\bigr)^2+E\sum_{y=0}^\infty \bbn_\nu(y) \mathop{\rm var}\nolimits \Bigl(\Lambda|Y=y\Bigr)
\\
&= E \sum_{y=0}^\infty  \Bigl((y+1)^2 \frac{\bbn_\nu^2(y+1)}{\bbn_\nu(y)}-2 (y+1)^2 \frac{\bbn_\nu(y+1)P(y+1)}{P(y)}
\\
&\qquad\qquad +\bbn_\nu(y) {\delta^G}^2(y)\Bigr)
1(\bbn_\nu(y)>0)+r(G,\delta^G)
\\
&= r(G,\delta^G) + E \sum_{y=0}^\infty  (y+1)^2 \frac{\bbn_\nu^2(y+1)}{\bbn_\nu(y)}1(\bbn_\nu(y)>0) - \nu E{\delta^G}^2(Y).
\end{split}
\end{equation*}

In the above we used the facts  that $\bbn_\nu(y+1)$ and $\bbn_\nu(y)$ are independent, and that if $X\sim Po(\theta)$ then $EX^2=\theta+\theta^2$.

In order to evaluate $ R(G,\hat{\delta}^G)-R(G,\delta^G)$ we need the following.
\begin{equation*}
E_\theta \frac{1(X>0)}X  = e^{-\theta} \sum_{i=1}^\infty \frac{\theta^i}{i!i},
\end{equation*}
hence \begin{equation}
\label{EX2}
\begin{split}
    E \frac{1(X>0)}X &= e^{-\theta}\summ i1\en\frac{\theta^i}{i!i}
    \\
    &= c e^{-\theta} \summ i1\en\frac{\theta^i}{(i+1)!}
    \\
    &= c   e^{-\theta} \theta^{-1}(e^{\theta} - 1 -\theta),
\end{split}
\end{equation}
where $c\in(1,2)$.

Also,
\begin{equation}
    \label{EX3}
    \begin{split}
       E \frac{1(X>0)}X-\frac1\theta
       &= e^{-\theta}\sum_{i=1}^\infty \frac{\theta^i}{i!i}-\frac1\theta
       \\
         &= e^{-\theta}\sum_{i=1}^\infty \frac{\theta^i}{(i+1)!} -\frac1\theta
         + e^{-\theta}\sum_{i=1}^\infty \frac{\theta^i}{(i+1)!i}
         \\
         &\le e^{-\theta}\sum_{i=1}^\infty \frac{\theta^i}{(i+1)!} -\frac1\theta
         + 3e^{-\theta}\sum_{i=1}^\infty \frac{\theta^i}{(i+2)!}
         \\
         &= \frac1\theta e^{-\theta} \Bigl(e^{\theta}-1-\theta\Bigr)-\frac1\theta+
         \frac3{\theta^2}e^{-\theta}\Bigl(e^{\theta}-1-\theta-\frac12\theta^2\Bigr)
         \\
         &= - \frac{1+\theta}\theta e^{-\theta} + \frac3{\theta^2}e^{-\theta}\Bigl(e^{\theta}-1-\theta-\frac12\theta^2\Bigr).
     \end{split}
\end{equation}

Now
\begin{equation*}
\begin{split}
r(G,\hat {\delta}^G) &=r(G,\delta^G) + \sum_{y=0}^\infty (y+1)^2E\frac{\nu P(y+1)}{\bbn_\nu(y)}1(\bbn_\nu(y)>0)
\\&+ \sum_{y=0}^\infty (y+1)^2 E\Bigl( \frac{\nu^2P^2(y+1)}{\bbn_\nu(y)}1(\bbn_\nu(y)>0) -\nu\frac{P^2(y+1)}{P(y)}\Bigr)
\\
&= r(G,\delta^G)+I+II, \quad\text{say}.
\end{split}
\end{equation*}
In the following   $c_1,\dots,c_5\in (a,b)$ are some constants for some universal constants $0<a<b<\en$. Now,
\begin{equation*}
\begin{split}
I=c_1 \sum_{y=0}^\infty (y+1)^2\frac{P(y+1)}{P(y)}\bigl(1-e^{-\nu P(y)}(1+\nu P(y))\bigr)
\end{split}
\end{equation*}
If $G^k$ has a compact support, then $\delta^{G^k}(y)$ is increasing and bounded by $\lambda_U \equiv \lambda_U^k<\lambda_2$, the upper support of $G^k$.  Using this observation and \eqref{EX2}, we obtain  for  $\nu$ large enough
\begin{equation*}
\begin{split}
I &= c_1 \lambda_U \sum_{y=0}^\infty (y+1) \bigl(1-e^{-\nu P(y)}(1+\nu P(y))\bigr)
\\
&=  c_1 \lambda_U \sum_{\nu P(y)>1/2 } (y+1) \bigl(1-e^{-\nu P(y)}(1+\nu P(y))\bigr)
\\
&\hspace{3em}+c_1 \lambda_U \sum_{\nu P(y)\le 1/2} (y+1) \bigl(1-e^{-\nu P(y)}(1+\nu P(y))\bigr)
\end{split}
\end{equation*}
Note that for $\theta>0$, $1-(1+\theta)e^{-\theta}$ is monotone increasing from 0  to 1:
\begin{equation*}
\begin{split}
I&= c_2 \lambda_U \sum_{\nu P(y)>1/2 } (y+1) + c_3  \sum_{\nu P(y)\le1/2} (y+1) \bigl(\nu P(y)\bigr)^2
\\
&= c_2\lambda_U\max\{(y+1)^2: P(y)>1/2\nu\}  + c_3 \lambda_U  \sum_{\nu P(y)\le1/2} (y+1) \bigl(\nu P(y)\bigr)^2
\end{split}
\end{equation*}
Now, for $z>2\lambda_U$, $\sum_{y\ge z}y^kP(y)\le 2 P(z)$, $k=0,1$, and $e^{-\lambda}\lambda^y/y!=\epsilon$ implies that $y \log|\log\epsilon|/|\log\epsilon|\to 1$ as $\epsilon\searrow0$ and $y \rightarrow \infty$ for any $\lambda_1\le\lambda\le\lambda_2$. Hence
\begin{equation*}
\begin{split}
I&= c_2 \lambda_U \bigl(\log \nu/\log\log\nu\bigr)^2 + c_3 \nu\sum_{\nu P(y)\le1/2} (y+1) P(y)
\\
&=c_4 \lambda_U  \bigl(\log \nu/\log\log\nu\bigr)^2 .
\end{split}
\end{equation*}

Bounding $II$ is similar, noting that there is $\gamma>0$ such that the RHS of \eqref{EX3} is negative for $\theta<\gamma$:
\begin{equation*}
\begin{split}
II&\leq  \sum_{y-0}^\infty (y+1)^2 \nu^2 P^2(y+1) E\Bigl( \frac{1(\bbn_\nu(y)>0)}{\bbn_\nu(y)}
- \frac1{\nu P(y)}\Bigr)
\\
&\leq 3\sum_{\nu P(y)>\gamma}(y+1)^2 \frac{P^2(y+1)}{P^2(y)}\Bigl(1-e^{-\nu P(y)}(1+\nu P(y)+\frac12\nu^2P^2(y))\Bigr)
\\
&\leq c_5 \sum_{\nu P(y)>\gamma}(y+1)^2 \frac{P^2(y+1)}{P^2(y)}
\\
&\leq c_6 \max\{y:\;\; \nu P(y)>\gamma\}
\\
&= c_6  \log \nu/\log\log\nu.
\end{split}
\end{equation*}

\end{proof}

\bigskip

{\bf Remarks:}
\begin{enumerate}
\item The asymptotics in the above theorem implies that in a non-sparse situation,
asymptotically there is
a room for only a negligible improvement on Robbins' classical estimator.
However, in light of Example 1 and our simulations,
the asymptotic presented in this section may be somewhat misleading.
This is since the above asymptotics often seems to `kick-in' only for  very large $n$ and
are thus irrelevant
for moderately large values of $n$, that appear in practice.


\item Nevertheless, our asymptotics suggests that there are limitations
and possible room for improvement of  Robbins' classical procedure in a triangular array setup of sparse problems
in which the risk may be of order $O((\log(n)/\log\log(n))^2))$, for arbitrarily large $n$.
\end{enumerate}

{\bf Acknowledgement.} We are grateful to the reviewers for their comments and suggestions
and in particular for Table 7 that was provided to us.

\vspace{3ex}
{\bf \Large References}

\begin{list}{}{\setlength{\itemindent}{-1em}\setlength{\itemsep}{0.5em}}

\item
Brown, L. D., Gans, N., Mandelbaum, A., Sakov, A., Shen, H.,
Zeltyn, S., and Zhao, L. H. (2005).
Statistical analysis of a telephone call center: a queing
science perspective. {\it Jour. Amer. Stat. Asoc.}
{\bf 100} 36-50.
\item
Brown, L.D., Cai, T., Zhang, R., Zhao, L., Zhou, H. (2010). The
root-unroot algorithm for density
estimation as implemented via wavelet block thresholding. {\it Probability and Related Fields}, {\bf 146}, 401-433.
\item
Brown, L.D. and Greenshtein, E. (2009). Non parametric
empirical Bayes and compound decision
approaches to estimation of high dimensional vector of normal
means. {\it Ann. Stat.} {\bf 37}, No 4, 1685-1704.
\item
Clevenson, L. and Zidek, J.V. (1975). Simultaneous estimation
of the mean of independent Poisson laws.
{\it Jour. Amer. Stat. Asoc.} {\bf 70} 698-705.
\item
Cohen, N., Greenshtein E., and Ritov, Y. (2012). Empirical Bayes in the presence of explanatory
variables. To appear in {\it Statistica Sinica.}
\item
Copas, J.B. (1969). Compound decisions  and empirical Bayes (with discussion). {\it JRSSB} {\bf 31} 397-425.
\item
Efron, B. (2003). Robbin, empirical Bayes, and microarrays (invited paper). {\it Ann.Stat.}
{\bf 31}, 364-378.
\item
Fourdrinier, D. and Robert, C. P. (1995). Intrinsic losses for empirical Bayes estimation:
A note on normal and Poisson cases. {\it Stat. and Prob. Letters} {\bf 23}, 35-44.
\item
Greenshtein, E. (1996). Comparison of sequential experiments.
{\it Ann. Stat.} {\bf 24}, No 1, 436-448.
\item
Greenshtein, E. and Ritov, Y. (2008). Asymptotic efficiency of
simple
decisions
for the compound decision problem. {The 3'rd
Lehmann Symposium. IMS Lecture
Notes Monograph Series}, J.Rojo, editor. 266-275.

\item  Hengartner, N. W. (1997). Adaptive Demixing in Poisson Mixture Models. \emph{Ann. of Statist.} \textbf{25}, 917--928.
\item
Jiang, W. and Zhang, C.-H. (2009). General maximum likelihood
empirical Bayes estimation of normal means. {\it Ann. Stat.} {\bf 37}, No 4, 1647-1684.
\item
Jiang, W. and Zhang, C.-H. (2010). Empirical Bayes in-season prediction of baseball batting average. {\it Borrowing Strength: Theory Powering Application-A festschrift for L.D. Brown}
J.O. Berger, T.T. Cai, I.M. Johnstone, eds. IMS collections  {\bf 6}, 263-273.
\item
Johnstone, I. (1984). Admissibility, difference equations and
recurrence in estimating a Poisson mean.
{\it Ann. Stat.} {\bf 12}, 1173-1198.
\item
Johnstone, I. and Lalley, S. (1984). On independent statistical
decision problems and products of diffusions.
{\it Z. fur Wahrsch.} {\bf 68}, 29-47.
\item
Koenker, R. and Mizera, I. (2012). Shape constraints , compound decisions and empirical Bayes rules. Manuscript.
\item
Kiefer, J. and Wolfowitz, J. (1956). Consistency of the maximum likelihood estimator in the presence of infinitely many incidental parameters. {\it Ann.Math.Stat.} {\bf 27 }, 887-906.
\item
Mammen, E. (1991). Estimating a smooth monotone Regression function. {\it Ann.Stat.} {\bf 19},
No. 2, 724-740.
\item
Maritz, J.S. (1969). Empirical Bayes estimation for the Poisson
distribution. {\it Biometrika} {\bf 56}, N0.2, 349-359.
\item
Park, J. (2011). Non parametric empirical Bayes estimator in simultaneous estimation
of Poisson means with application to mass spectrometry data. {\it Journal of Nonparametric Statistics}.
\item
R Development Core Team (2008). R: A language and environment for
  statistical computing. R Foundation for Statistical Computing,
  Vienna, Austria. ISBN 3-900051-07-0, URL http://www.R-project.org.
\item
Robbins, H. (1951). Asymptotically subminimax solutions  of
compound decision problems. {\it Proc. Second Berkeley Symp.}
131-148.
\item
Robbins, H. (1955). An Empirical Bayes approach to statistics.
{\it  Proc. Third Berkeley Symp.} 157-164.
\item
Robbins, H. (1964). The empirical Bayes approach to statistical
decision problems. {\it Ann.Math.Stat.}
{\bf 35}, 1-20.
\item
Robertson, T., Wright, F. T. and Dykstra, R. L. (1988).
\emph{Order Restricted Statistical Inference}. Wiley,
New York.
\item
Zhang, C.-H.(2003). Compound decision theory and empirical
Bayes methods.(invited paper). {\it Ann. Stat.} {\bf 31}
379-390.
\item
\end{list}

\end{document}